\documentclass[12pt, reqno]{amsart}

\usepackage{amsthm, amssymb}
\usepackage{mathrsfs}

\usepackage{enumerate}

\usepackage[colorlinks]{hyperref}

\allowdisplaybreaks

\theoremstyle{plain}
\newtheorem{thm}{Theorem}[section]
\newtheorem{defn}[thm]{Definition}
\newtheorem{lem}[thm]{Lemma}
\newtheorem{prop}[thm]{Proposition}
\newtheorem{cor}[thm]{Corollary}

\title[Paley inequality for the Weyl transform]{Paley inequality for the Weyl transform and its applications}
\author{Ritika Singhal}
\address{Ritika Singhal \newline \hspace*{0.25cm} 
        Department of Mathematics\newline \hspace*{0.25cm}
	Indian Institute of Technology Delhi\newline \hspace*{0.25cm}
	Delhi - 110016\newline \hspace*{0.25cm}
	India.}
\email{ritikasinghal1120@gmail.com}
\author{N. Shravan Kumar}
\address{N. Shravan Kumar \newline \hspace*{0.25cm}
        Department of Mathematics\newline \hspace*{0.25cm}
	Indian Institute of Technology Delhi\newline \hspace*{0.25cm}
	Delhi - 110016\newline \hspace*{0.25cm}
	India.}
\email{shravankumar.nageswaran@gmail.com}

\begin{document}
	
\begin{abstract}
    In this paper, we prove several versions of the classical Paley inequality for the Weyl transform. As an application, we  discuss  $L^p$-$L^q$ boundedness of the Weyl multipliers and prove a version of the H\"ormander's multiplier theorem. We also prove Hardy-Littlewood inequality. Finally, we study vector-valued versions of these inequalities. In particular, we consider the inequalities of Paley, Hausdorff-Young, and Hardy-Littlewood and their relations.\end{abstract}

	

\maketitle

\section{Introduction}
Let $G$ be a locally compact abelian group. A classical result of the Fourier analysis is the Hausdorff-Young inequality which states the following: if $1 \leq p \leq 2$, then the Fourier transform maps $L^p(G)$ into $L^{p'}(\widehat{G})$, where $\widehat{G}$ is the dual group of $G$ and $p'$ is the conjugate index of $p$. Later, Paley \cite{Paley} extended this result to Lorentz spaces and showed that for $G=\mathbb{T}$, if $f \in L^p(\mathbb{T})$, then $\hat{f} \in l^{p',p}(\mathbb{Z})$. The case $G=\mathbb{R}$ is due to H\"ormander \cite{Hor}.

The Weyl transform, defined by Hermann Weyl in \cite{Weyl}, is a pseudo-differential operator associated to a measurable function 
on $\mathbb{R}^n \times \mathbb{R}^n$. A large number of both mathematicians and physicists have studied the properties of the Weyl transform and its applications to quantum mechanics and partial differential equations. The Fourier transform and the Fourier inversion formula serves as the foundation for the creation of the Weyl transform. Therefore, all the classical properties of the Fourier transform namely the Reimann-Lebesgue lemma, Plancheral theorem, and Hausdorff-Young inequality, works for the Weyl transform as well. A natural question that arises here is whether Paley's extension of the Hausdorff-Young inequality possible for the Weyl transform. This paper answers this question in the affirmative sense. In fact, the main aim of this paper is to study Paley's inequality and its variants for the Weyl transform associated to locally compact abelian groups. More precisely, we  derived  the Paley inequality for the Weyl transform and the Inverse Weyl transform  in Sections \ref{HYPI1} and  \ref{IWT}, respectively, and proved the following version of the Paley inequality  in Theorem  \ref{PIWT}.

\begin{thm}[Paley inequality]
		Consider a positive function $\phi \in l^{1, \infty}(\mathbb N)$. Then for $ 1 <p \leq 2$ and  $f\in L^p(G\times\widehat{G})$, we have $$\left( \sum S_n(W(f))^p\phi(n)^{2-p}\right)^\frac{1}{p} \lesssim \|\phi\|_{l^{1,\infty}(\mathbb N)}^{\frac{2-p}{p}}\|f\|_p.$$
 In particular, for  $f \in L^p( G \times \widehat{G}),$ we have $ W(f) \in \mathcal{B}_{p',p}(L^2(G)).$
	\end{thm}
Interpolating  the Hausdorff-Young inequality with  Paley-inequality, one can obtain the  Hausdorff-Young-Paley
inequality for the Weyl transform which is discussed in Theorem \ref{HYPI}.

The study of Fourier multipliers on $L^p$ spaces is one of the classical topics of harmonic analysis. One of the celebrated results of H\"ormander gives sufficient conditions for a symbol to be a $L^p$-$L^q$ Fourier multiplier. Very recently, Ruzhansky with his coauthors  studied the H\"ormander's Theorem for  Lie groups\cite{R1,ANR}, homogeneous manifolds\cite{R4}, locally compact groups\cite{AR} and compact hypergroups\cite{R2}. In \cite{Mau}, Mauceri introduced the concept of Weyl multipliers and proved a version of H\"ormander's theorem for the Weyl transform on $\mathbb{R}^{2n}$. The assumptions include some regularity conditions like the commutator of the operator with annihilation and creation operator should be a bounded operator. In this paper, as an application of the Paley inequality, we also prove a version of the H\"ormander's theorem for the Weyl transform on locally compact abelian groups without using any regularity assumptions. Thus, our conditions are in terms of the  singular value sequence associated to the operator and we have proved the following result: 
	\begin{thm}[H\"{o}rmander's theorem] 
	Let $G$ be a locally compact abelian group and let $1<p\leq 2\leq q<\infty.$ For $M \in \mathcal{B}(L^2(G)),$ 
  consider the operator $C_M$ defined on $L^1 \cap L^2(G\times \widehat{G} )$ as $W(C_Mf)=MW(f). $
  If $\{S_n(M)\} \in  l^{r,\infty}(\mathbb N)$  where $\frac{1}{r}=\frac{1}{p}-\frac{1}{q}$ then $M \in \mathscr{M}_{p,q}$ and 
  $$\|C_Mf\|_{L^q(G\times\widehat{G})} \lesssim \underset{s>0}{\sup}\ \left[s \underset{S_n(M)>s}{\underset{n\in\mathbb{N}}{\sum}}1  \right]^{\frac{1}{p}-\frac{1}{q}}\|f\|_p.$$
	\end{thm}

 Yet another classical inequality of Fourier analysis is the Hardy-Littlewood inequality. Established by Hardy and Littlewood for the torus $\mathbb{T}$ in (\cite{HL}), it  was later extended by 
Hewitt and Ross \cite{HR} for compact abelian groups. In 2016, this inequality was further extended to compact Lie groups in \cite{ANR} and to locally compact separable unimodular groups by Kosaki \cite{Kos}. In \cite{R2}, the same was proved as an application of Paley inequality. In Section \ref{HLI}, as another application of the Paley inequality, we prove an analogue of the Hardy-Littlewood inequality for the Weyl transform on locally compact abelian groups.
 
 Throughout the past few decades, the study of vector-valued functions has become increasingly popular. A majority of the classical problems in the theory of functions may be investigated in a vector-valued environment. Such a study provides new ideas for understanding various mathematical problems. As a result, several researchers have  emphasized its importance and extended the classical results to the vector-valued 
 setting- Hausdorff-Young inequality \cite{Bour1, GKKT, GMP, jaak, HHL, VKSK}, Paley inequality \cite{GKK}, Hardy's inequality \cite{B1}, singular integrals \cite{Bour2}, Fourier multipliers \cite{JSRR, GJR, SSV, RV} and Weyl transform \cite{RSSK}.
 
The vector-valued notion of Paley inequalities was first studied by Garcia et. al in \cite{GKK, GKKT}. Recently, in \cite{DM}, the authors study vector-valued versions of the Hardy-Littlewood inequality and the Paley inequality. In Sections \ref{HLI} and \ref{PTI}, we introduce the concept of Weyl-Paley type/cotype and Weyl-HL type/cotype for a Banach space $X$ and prove an analogue of Hausdorff-Young inequality. Finally, we examine the relationships between the Weyl type introduced in \cite{RSK}, the Paley type, and the HL type (as well as, their cotype counterparts) with each other. 

We shall begin with some preliminaries that are needed in the sequel.

\section{Preliminaries}
Let $G$ be a locally compact abelian group with $\widehat{G}$ as its dual group. As usual, for $1\leq p \leq \infty$, we shall denote by $L^p(G)$, the usual classical $L^p$-spaces on $G$ w.r.t. the Haar measure on the group $G$. We shall denote by $\mathcal{B}(L^2(G))$,  the Banach algebra of all bounded operators  on $L^2(G)$ and $C_c(G)$ will denote all compactly supported functions on the group $G$.

The \textit{Weyl transform}, denoted  $W,$ is defined as a $\mathcal{B}\left(L^{2}(G)\right)$-valued integral on $C_{c}(G \times \widehat{G})$ given by $$W(f) \varphi(y)=\int_{G \times \widehat{G}}f(x, \chi) \rho((x, \chi))(\varphi)(y) dx  d\chi, ~f \in C_c(G \times \widehat{G})$$ where $\rho(x,\chi)$ is the \textit{Schr\"odinger representation} of $G \times \widehat{G}$ on $L^2(G)$, defined as $$\rho_{G}((x, \chi))(\varphi)(y)= \chi(y) \varphi(x y), \varphi \in L^{2}(G).$$

Let $\mathcal{H}$ be a  Hilbert space and $1\leq p<\infty.$ If  $T:\mathcal{H}\rightarrow \mathcal{H}$ is a compact operator, then it admits an orthonormal representation $$T=\underset{n\in\mathbb{N}}{\sum}S_n(T)\langle.,e_n\rangle \sigma_n,$$ where $\{e_n\}$ and $\{\sigma_n\}$ are orthonormal sequences in $\mathcal{H}$ and $S_n(T)$ denotes the $n^{th}$ singular value of $T$. The \textit{$p^{th}$-Schatten-von Neumann class}, denoted $\mathcal{B}_p(\mathcal{H}),$ consists of all compact operators, $T:\mathcal{H}\rightarrow\mathcal{H}$ with $\{S_n(T)\}\in\ell^p.$
 For $T\in\mathcal{B}_p(\mathcal{H}),$ define $$\|T\|_{\mathcal{B}_p(\mathcal{H})}:=\left(\underset{n\in\mathbb{N}}{\sum}|S_n(T)|^p\right)^{1/p}= \|\{S_n(T)\}\|_{l^p}.$$ The space $\mathcal{B}_p(\mathcal{H})$ with the above norm becomes a Banach space. We shall denote by $\mathcal{B}_\infty(\mathcal{H})$, the space of all compact operators on $\mathcal{H}.$ By Riesz-Thorin theorem on complex interpolation, for $0<\theta<1,$ we have, $[\mathcal{B}_1(\mathcal{H}),\mathcal{B}(\mathcal{H})]_\theta\cong\mathcal{B}_{1/\theta}(\mathcal{H}).$ For $1 \leq p \leq \infty$, $p'$ denoted the conjugate index of $p$ such that $1 / p+1 / p^{\prime}=1$.
 \begin{thm}
     Let $1 \leq p \leq 2$, then the Weyl transform is a continuous mapping of functions $f \in L^p(G \times \widehat{G} )$ to operators $W(f) \in \mathcal{B}_{p'}(L^2(G))$, i.e., there exists a constant $C>0$ such that
$$
 \left\|W(f) \right\|_{\mathcal{B}_{p'}(L^2(G)} \leqslant C\|f\|_{L^p(G \times \widehat{G})},  f \in L^p(G \times \widehat{G}).
$$
In fact for $p=2$, the Weyl transform is a unitary map between $L^2(G\times\widehat{G})$ and $\mathcal{B}_2(L^2(G)).$
\end{thm}

Also, by  duality, we can conclude that for $1 \leq p \leq 2$, if  
$W(f)$ belongs to $\mathcal{B}_{p}(L^2(G ))$ for some measurable function $f$, then $f$ belongs to $L^{p'}(G \times \widehat{G})$, and there exists $C>0$ such that
$$
 \|f\|_{L^{p'}(G \times \widehat{G})} \leq  C \left \|W(f) \right\|_{\mathcal{B}_{p}(L^2(G)}.
$$
For more on Weyl transform and Schatten-class operators, see \cite{Fol1, Wong, S}.

\begin{defn}
    Let $1 \leq p \leq \infty$. A bounded operator $M \in \mathcal{B}(L^2(G))$  is said to be a (left) Weyl multiplier of $L^p(
    G \times \widehat{G})$ if the operator $C_M $ defined on $f \in L^1\cap L^2(G \times \widehat{G})$ by $W(C_Mf)=MW(f)$ extends to a bounded operator on $L^p(G \times \widehat{G})$.
\end{defn}

Let $1 \leq p,q \leq \infty$. If $M \in \mathcal{B}(L^2(G))$ is a Weyl multiplier of $L^p(
    G \times \widehat{G})$, we say that $ M \in \mathscr{M}_{p,q}  $ if the map $C_M$ defined above extends to a bounded map from  $L^p(G \times \widehat{G})$ to $L^q(G \times \widehat{G})$.
In the case, when $p=q$,
$\mathscr{M}_{p,p}=\mathscr{M}_{p}$. By the Plancheral formula, $\mathscr{M}_2$ is the algebra $\mathcal{B}(L^2(G))$ and $\mathscr{M}_1$ coincides with the algebra of Weyl transform  of the space of finite Borel measures on $G \times \widehat{G}$. See \cite{RSK}.

Let $X$ be a Banach space and let $(\Omega, \mathcal{A})$ be a measure space with a $\sigma$-finite positive measure $\mu$. Consider a weight function $\omega: \Omega \to (0, \infty)$ which is integrable on sets of finite measure from $\mathcal{A}.$ For $p \in [1,\infty)$, we
define the Bochner spaces $L^p(\Omega,\omega ; X)$ formed by all (equivalence classes of) strongly
$\mu$-measurable functions $f : \Omega \to X$ having a finite norm
$$\|f\|_{L^p(\Omega, \omega;X)}=\left(\int_\Omega \|f(x)\|^p_{X} \omega(x) d \mu(x)\right)^{\frac{1}{p}}.$$ When $p = \infty, L^\infty(\Omega, \omega; X)$ = $L^\infty(\Omega; X)$ denote the functions which are essentially bounded and $$\|f\|_{L^\infty(\Omega, \omega, X)}= \text{ess.} \sup_{x \in \Omega}\|f(x)\|_{X}.$$ Also, $L^p(\Omega; X)$  denotes the  special case when $\omega \equiv 1$ and also $L^p(\Omega, \omega) = L ^p(\Omega, \omega; \mathbb{C}).$

\begin{defn}
    Let $f: \Omega \to X$. The decreasing rearrangement of $f$ is the function $f^*$ defined on $(0, \mu(\Omega))$ by 
    $$f^*(t)= \inf \{s>0: d_f(s) \leq t\}$$
     where $d_f(s )= \mu (\{x \in \Omega: \|f(x)\|_{X}> s\})$, the distribution function of $f$.
\end{defn}
\begin{defn}(Lorentz spaces)
    For $1 \leq p,q \leq \infty$, the Lorentz space $L^{ p,q}(\Omega; X)$ is the class of all $\mu$-strongly measurable
functions $f : \Omega \to X$ such that
$$\|f\|_{L^{p,q}(\Omega;X)}= \begin{cases} \left( \int_{0 }^{\mu(\Omega)}(t^{\frac{1}{p}}f^*(t))^q \  \frac{dt}{t}\right)^\frac{1}{q}< \infty & \text{if} \ q< \infty  \\ \underset{t>0}{\sup}\ t^{\frac{1}{p}}f^*(t) < \infty & \text{if} \  q= \infty. \end{cases} $$

\end{defn}
In the case when $\Omega $ is countable with discrete measure, we denote the above space by $l^{p,q}(\Omega;X)$. For $X= \mathbb{C}$, Lorentz spaces will be denoted by $L^{p,q}(\Omega).$ 
For all $0 < p,q \leq \infty$, the
spaces $L^{p,q}(\Omega, X)$ are complete with respect to their quasinorm and they are therefore quasi-Banach spaces. For Schatten classes, by replacing $l^p$-norm of the singular values by the Lorentz spaces $l^{p,q}$ quasi-norm, we get the non commutative Lorentz-spaces $\mathcal{B}_{p,q}(\mathcal{H})$ defined as the space of all compact operators $T \in \mathcal{B}(\mathcal{H})$ such that
$$\|T\|_{p,q}=\|\{S_n(T)\}\|_{l^{p,q}}= \begin{cases}
    \left(  \underset{n}{\sum} (n^{1/p-1/q} S_n(T))^q  \right )^{1/q} , & 1 \leq q < \infty,\\ \underset{n}{\sup} \  n^{1/p}S_n(T), & q= \infty
\end{cases} $$ is finite. The operators that map $L^p(\Omega)$ to $L^{q,\infty}(\Omega)$ are called weak type $(p,q)$.

 Let $X$ and $Y$ be Banach spaces. Let $ X\overset{\vee}{\otimes}Y$ and $X\overset{\wedge}{\otimes}Y$  denote the injective and projective tensor product   of $X$ and $Y$ respectively.
 If $\mathcal{H}$ is a Hilbert space, we shall denote by $\mathcal{B}_1[\mathcal{H};X]$ and $\mathcal{B}_\infty[\mathcal{H};X]$ the spaces $\mathcal{B}_1\overset{\wedge}{\otimes}X$ and $\mathcal{B}_\infty\overset{\vee}{\otimes}X,$ respectively. For $1<p,r<\infty,$ we define $ \mathcal{B}_p[\mathcal{H};X]$ and $\mathcal{B}_{p,r}[\mathcal{H};X]$ as follows:
 \begin{align*}
     \mathcal{B}_p[\mathcal{H};X]:=& \left[\mathcal{B}_1[\mathcal{H};X],\mathcal{B}_\infty[\mathcal{H};X]\right]_{1/p} \mbox{(in the sense of complex interpolation)} \\ \mathcal{B}_{p,r}[\mathcal{H};X]:=& \left[\mathcal{B}_1[\mathcal{H};X],\mathcal{B}_\infty[\mathcal{H};X]\right]_{1/p,r} \mbox{(in the sense of real interpolation)}.
 \end{align*}
 A set $Y \subseteq  X^*$ is norming for $X$ if $\sup_{f \in Y  \setminus \{0\} } \frac{|f(x)|}{\|f\|_{X^*}}=\|x\|_X.$

\begin{lem}
    Let $(\Omega, \mathcal{A}, \mu)$ be a $\sigma$-finite measure space and $\omega: \Omega \to (0, \infty)$ be measurable. Let $p \in (1, \infty)$ and let $Y \subseteq X^*$ be a normed closed subspace of $X^*$ which is norming for
$X$. Then $L^{p'}(\Omega,\omega^{\frac{1}{p-1}}; Y )$ is norming for $L^p(\Omega, \omega; X)$ with respect to the duality pairing
 $$\langle f,g \rangle =\int_\omega \langle f(x), g(x)\rangle d\mu(x).$$
Moreover, the subspace of $Y$-valued simple functions in $L^{p'}(\Omega, \omega^{\frac{1}{p-1}}; Y )$ is also norming for $L^p(\Omega, \omega; X)$.
\end{lem}
 The proof of the above lemma can be found in \cite[Proposition 1.3.1]{MR3617205}  for the unweighted case and the weighted case can then be generalised. 
 \begin{lem}[\cite{DM},Lemma 2.5]
      Let $(\Omega, \mathcal{A}, \mu)$ be a $\sigma$-finite non atomic measure space. Let $p,q \in (1, \infty)$ and let $Y \subseteq X^*$ be a normed closed subspace of $X^*$ which is norming for $X$. Then 
      $$\|f\|_{L^{p,q}(\Omega ;X)} \cong \sup \left \{  \left|  \int_\Omega \langle f(x),g(x) \rangle  d \mu(x) \right|: \|g\|_{L^{p', q'}(\Omega;Y)} \leq 1 \right \}.$$
 \end{lem}
 
 \begin{lem}[\cite{Pis}]
     Let $1<p< \infty$. Then, we have 
         $$\mathcal{B}_p[\mathcal{H};X]^*\cong \mathcal{B}_{p'}[\mathcal{H};X^*].$$
 \end{lem}
 \begin{defn}
     		We say that a Banach space $X$ has Weyl type $p$ if there exists a constant $C>0$ such that
       $$\|W(f)\|_{\mathcal{B}_{p^\prime}[L^2(G);X]} \leq C \|f\|_{L^p(G\times\widehat{G},X)}.$$
       Similarly, 
        a Banach space $X$ has Weyl cotype $q$ if there exist a constant $C>0$ such that
       $$\|f\|_{L^q(G\times\widehat{G},X)}  \leq C \|W(f)\|_{\mathcal{B}_{q^\prime}[L^2(G);X]} .$$

    \end{defn}
For more on Lorentz spaces and interpolation theorems, one can refer to \cite{BS} or \cite{Gra}. The following is an adaptation of the classical Marcinkiewicz interpolation theorem for Schatten class operators. Since we couldn't find the result anywhere, we are providing a brief proof of it. For more general versions see \cite{JX, CR}.
    \begin{thm}\label{NCommCommMIT}
        Let $(X,\mu)$ be a measure space and let $\mathcal{H}$ be a separable Hilbert space. For each $0<p_0<p_1<\infty,$ let $\varphi:\mathcal{B}_{p_0}(\mathcal{H})+\mathcal{B}_{p_1}(\mathcal{H})\rightarrow L^0(X,\mu)$ be a sublinear operator. Suppose that there exist constants $C_0$ and $C_1$ such that 

        \begin{equation} \label{eqn 19}
              \|\varphi(T)\|_{L^{p_i,\infty}}\leq  C_i \|T\|_{\mathcal{B}_{p_i}(\mathcal{H})}\ \forall\ T\in\mathcal{B}_{p_i}(\mathcal{H}),  i=0,1 
        \end{equation}
       
        Then $\forall$ $p_0<p<p_1$ and $\forall$ $T\in\mathcal{B}_p(\mathcal{H}),$ there exists $C>0$ such that $$\|\varphi(T)\|_{L^p}\leq C\|T\|_{\mathcal{B}_p(\mathcal{H})}.$$
    \end{thm} 
    \begin{proof}
        Fix $T\in \mathcal B_p(\mathcal{H})$ and $\alpha>0.$ Since $T$ is a compact operator, $\exists$ orthonormal sequences $\{\sigma_n\}$ and $\{\theta_n\}$ for $\mathcal{H}$ such that $$T=\underset{n\in\mathbb{N}}{\sum}\lambda_n(\sigma_n\otimes\theta_n).$$ For each $n\in\mathbb{N},$ let $$\lambda_{0,n}^\alpha=
        \left\{ \begin{array}{cc} 
        \lambda_n & \mbox{if } |\lambda_n| >\delta\alpha \\
        0 & \mbox{otherwise} \end{array} \right. 
        \mbox{ and } \lambda_{1,n}^\alpha=
        \left\{ \begin{array}{cc} 
        \lambda_n   & \mbox{if } |\lambda_n| \leq\delta\alpha \\ 
        0 & \mbox{otherwise} \end{array} \right.$$ for some fixed $\delta \text{and} \alpha >0.$ Let $T_0^\alpha=\underset{n\in\mathbb{N}}{\sum}\lambda_{0,n}^\alpha (\sigma_n\otimes\theta_n)$ and $T_1^\alpha=\underset{n\in\mathbb{N}}{\sum}\lambda_{1,n}^\alpha (\sigma_n\otimes\theta_n).$ Then $T=T_0^\alpha+T_1^\alpha$ and $|\varphi(T)|\leq |\varphi(T_0^\alpha)| + |\varphi(T_1^\alpha)|.$ Now we omit the remaining part of the proof as it goes exactly as in the classical case. See \cite{Gra}.
    \end{proof}
	
\section{Hausdorff-Young Paley inequality} \label{HYPI1}

    In this section, we prove the Hausdorff-Young Paley inequality for the Weyl transform.
    We shall begin this section by proving the Paley inequality. The corresponding analogue of this for the Fourier transform on $\mathbb{R}$ can be found in \cite{Hor}.
	\begin{thm}\label{PIWT}[Paley inequality]
		Consider a positive function $\phi \in l^{1, \infty}(\mathbb N)$. Then for $ 1 <p \leq 2$ and  $f\in L^p(G\times\widehat{G})$, we have $$\left( \sum S_n(W(f))^p\phi(n)^{2-p}\right)^\frac{1}{p} \lesssim \|\phi\|_{l^{1,\infty}(\mathbb N)}^{\frac{2-p}{p}}\|f\|_p.$$
	\end{thm}
	\begin{proof}
		Consider the measure $\nu $ on $\mathbb{N}$ given by 
		\begin{equation}\label{distribution}
			\nu(\{n\}):=\phi^2(n).
		\end{equation}
		For $1 < p \leq \infty$, we let $l^p(\mathbb{N},\nu)$ denote the space of all complex-valued sequences $x=(x_n)_{n\in\mathbb{N}}$ such that $\|x\|^p_p= \underset{n\in\mathbb{N}}{\sum}|x_n|^p\phi^2(n)<\infty.$ We now claim that if $f\in L^p(G\times\widehat{G}),$ then $\left\{ \frac{S_n(W(f))}{\phi(n)} \right\}_{n \in \mathbb{N}}\in\ell^p (\mathbb{N},\nu).$ We will denote this correspondence by $T$ and show that $T$ is a bounded map. Our strategy here is to make use of the classical Marcinkiewicz Interpolation theorem \cite[Theorem 1.3.2]{Gra}. To do this we first show that $T$ is both weak type $(2,2)$ and $(1,1).$
  
  The distribution function, in this case, is given by  $$d_{T(f)}(y)=\nu(\left\{ n\in\mathbb{N}:|T(f)(n)|>y \right\}).$$ 
		To show that $T$ is of weak type $(1,1),$ we prove that $$\|T(f)\|_{1,\infty} \lesssim \|\phi\|_{l^{1,\infty}(\mathbb N)} \|f\|_1.$$ Observe that $$S_n(W(f))\leq \underset{n\in\mathbb{N}}{\sup}  S_n(W(f)) \leq \|W(f)\|\leq \|f\|_1,$$ and therefore, for $y<\frac{S_n(W(f))}{\phi(n)} \leq \frac{\|f\|_1}{\phi(n)},$ we have $$\nu\left(\left\{n\in\mathbb{N}:\frac{S_n(W(f))}{\phi(n)}> y \right\} \right) \leq \nu\left(\left\{n\in\mathbb{N}:\frac{\|f\|_1}{\phi(n)}>y\right\}\right).$$ Hence $\underset{y<|T(f)(n)|}{\underset{n\in\mathbb{N}}{\sum}}\ \phi^2(n) \leq \underset{y<\frac{\|f\|_1}{\phi(n)}}{\underset{n\in\mathbb{N}}{\sum}}\ \phi^2(n).$
		
		Now, let $w=\frac{\|f\|_1}{y}.$ Then 
		\begin{eqnarray*}
			\underset{\phi(n)<w}{\underset{n\in\mathbb{N}}{\sum}} \phi^2(n) &=& \underset{\phi(n)<w}{\underset{n\in\mathbb{N}}{\sum}}\ \int_0^{\phi^2(n)} d\tau = \int_0^{w^2}\ d\tau\ \underset{\sqrt{\tau}<\phi(n) < w}{\underset{n\in\mathbb{N}}{\sum}} 1 \\ &=& \int_0^w 2s\ ds\ \underset{s<\phi(n)<w}{\underset{n\in\mathbb{N}}{\sum}}\ 1 \leq \int_0^w 2\left( s\underset{s<\phi(n)}{\underset{n\in\mathbb{N}}{\sum}}\ 1\right)\ ds \\ &\leq& \int_0^w\ 2 \|\phi\|_{ l^{1,\infty}(\mathbb N)}\ ds = 2w \|\phi\|_{ l^{1,\infty}(\mathbb N)} = \frac{2\|\phi\|_{ l^{1,\infty}(\mathbb N)}}{y}\|f\|_1.
		\end{eqnarray*}
		Thus, for $y>0$, we have $$yd_{Tf}(y) = y \underset{y<|T(f)(n)|}{\underset{n\in\mathbb{N}}{\sum}}\ \phi^2(n) \lesssim \|\phi\|_{ l^{1,\infty}(\mathbb N)}\|f\|_1$$
		 Also, by using the Plancheral theorem for Weyl transform, it can be seen that $T$ maps $L^2(G \times \widehat{G})$ continuously to $l^2(\mathbb{N},\nu)$ since 
  	\begin{eqnarray*}
			\underset{n\in\mathbb{N}}{\sum} |T(f)(n)|^2\phi^2(n) =& \underset{n\in\mathbb{N}}{\sum} |S_n(W(f))|^2 = \|W(f)\|^2_{\mathcal{B}_2(L^2(G))} = \|f\|^2_2. 
		\end{eqnarray*}
	This shows that $T$ is weak type $(2,2)$.
 
		Finally, using Marcinkiewicz interpolation theorem, it follows that $\|T(f)\|_p\lesssim \|\phi\|_{ l^{1,\infty}(\mathbb N)}^{\left( \frac{2-p}{p} \right)}\|f\|_p$ or $$\left(\underset{n\in\mathbb{N}}{\sum}\ S_n(W(f))^p\phi(n)^{2-p}\right)^{1/p}\lesssim \|\phi\|_{ l^{1,\infty}(\mathbb N)}^{\frac{2-p}{p}}\|f\|_p.$$ Hence the proof.
	\end{proof}
 \begin{cor} \label{cor}
      For $ 1 <p \leq 2$ and  $f\in L^p(G\times\widehat{G})$, we have $W(f) \in \mathcal{B}_{p',p}(L^2(G))$ and there exists $C>0$ such that
      $$\|W(f)\|_{\mathcal{B}_{p',p}(L^2(G))} \leq C \|f\|_p.$$
 \end{cor}
	By interpolating the  Paley inequality and the Hausdorff-Young inequality for the Weyl transform, we obtain the following Hausdorff-Young Paley inequality.
 
	\begin{thm}[Hausdorff-Young Paley inequality] \label{HYPI}
		Let $G$ be a locally compact abelian group and let $1<p\leq b\leq p^\prime<\infty.$ If  $\phi \in  l^{1,\infty}(\mathbb N)$, then for all $f\in L^p(G\times\widehat{G}),$ we have $$\left( \underset{n\in\mathbb{N}}{\sum}\left(S_n(W(f)) \phi(n)^{\frac{1}{b} -\frac{1}{p^\prime}} \right)^b\right)^\frac{1}{b}\lesssim \|\phi\|_{ l^{1,\infty}(\mathbb N)}^{\frac{1}{b} - \frac{1}{p^\prime}}\|f\|_p.$$
	\end{thm}
	We now prove, as an application, the Weyl transform analogue of the H\"{o}rmander's theorem.
	\begin{thm}[H\"{o}rmander's theorem] \label{Hor}
		Let $G$ be a locally compact abelian group and let $1<p\leq 2\leq q<\infty.$ For $M \in \mathcal{B}(L^2(G)),$ 
  consider the operator $C_M$ defined on $L^1 \cap L^2(G\times \widehat{G} )$ as $W(C_Mf)=MW(f). $
  If $\{S_n(M)\} \in  l^{r,\infty}(\mathbb N)$  where $\frac{1}{r}=\frac{1}{p}-\frac{1}{q}$ then $M \in \mathscr{M}_{p,q}$ and 
  $$\|C_Mf\|_{L^q(G\times\widehat{G})} \lesssim \underset{s>0}{\sup}\ \left[s \underset{S_n(M)>s}{\underset{n\in\mathbb{N}}{\sum}}1  \right]^{\frac{1}{p}-\frac{1}{q}}\|f\|_p.$$
	\end{thm}
 
 \begin{proof}
		Let $p\leq q^\prime.$ Then for $f\in C_c(G\times\widehat{G}),$ 
		\begin{eqnarray*}
			\|C_Mf\|_q \leq \|W(C_Mf)\|_{\mathcal{B}_{q'}(L^2(G))} &= & \left( \underset{n\in\mathbb{N}}{\sum} S_n(W(C_Mf))^{q^\prime}  \right)^{\frac{1}{q^\prime}} \\ &=& \left( \underset{n\in\mathbb{N}}{\sum} S_n(MW(f))^{q^\prime} \right)^{\frac{1}{q^\prime}}.
		\end{eqnarray*}
  Since $S_{n+m+1}(T) \leq S_{n+1}(T)+S_{m+1}(T)$ for any compact operator $T$, we have
  $$ \left( \underset{n\in\mathbb{N}}{\sum} S_n(MW(f))^{q^\prime} \right)^{\frac{1}{q^\prime}} \\ \lesssim  \left( \underset{n\in\mathbb{N}}{\sum} \left(S_n(M)S_n(W(f))\right)^{q^\prime} \right)^{\frac{1}{q^\prime}}.$$
		Let $b=q'$ and $\phi(n)=S_n(M)^r.$ Hence using Theorem \ref{HYPI}, we have 
		
		\begin{eqnarray*}
			\left( \underset{n\in\mathbb{N}}{\sum} \left(S_n(M)S_n(W(f))\right)^{q^\prime} \right)^{\frac{1}{q^\prime}} \lesssim  \left( \underset{s>0}{\sup}  \ s \underset{\underset{S_n(M)^r > s}{n \in \mathbb{N}}}{\sum} 1\right)^\frac{1}{r} \|f\|_p.
		\end{eqnarray*}

  Also, notice that
		$$
	\left( \underset{s>0}{\sup} \  s \underset{\underset{S_n(M)^r > s}{n \in \mathbb{N}}}{\sum} 1\right)^\frac{1}{r} =
 \left(\underset{s>0}{\sup} \  s^r \underset{\underset{S_n(M) > s}{n \in \mathbb{N}}}{\sum} 1\right)^\frac{1}{r} = \underset{s>0}{\sup} \  s \left( \underset{\underset{S_n(M) > s}{n \in \mathbb{N}}}{\sum} 1\right)^\frac{1}{r} $$

 Hence, we get 
 $$\|C_Mf\|_q \lesssim \underset{s>0}{\sup} \  s \left( \underset{\underset{S_n(M) > s}{n \in \mathbb{N}}}{\sum} 1\right)^{\frac{1}{p}-\frac{1}{q} } \|f\|_p$$
 as required.

	Now, by using the duality , we know that $$\|C_M\|_{L^p(G \times \widehat{G}) \to L^q(G \times \widehat{G})}=\|C_M^*\|_{L^{q'}(G \times \widehat{G}) \to L^{p'}(G \times \widehat{G})}.$$ Hence for the case when $q'<(p')'=p$, one can work with $C_M^*$ whose associated operator will be $M^*$ and $\|C_M^*\|=\|C_M\|$.
	\end{proof}
	
	\section{Hausdorff-Young inequality for the inverse Weyl transform} \label{IWT}
            In this section, we prove the Paley inequality for the inverse Weyl transform. Our approach here is the same as what we did for the case of the Weyl transform. Finally, we also prove a version of the H\"{o}rmander's theorem for the inverse Weyl transform which will later lead to Hardy-Littlewood inequality.

	\begin{thm} \label{Paley}

		 Consider a positive function $\psi \in L^{1, \infty}(G\times\widehat{G})$. Let 
		\begin{equation}\label{Paley-Eqn} 
			M_\psi:= \|\psi\|_{L^{1, \infty}(G\times\widehat{G})}=\sup_{s>0} \ s \underset{ \underset{|\psi(x,\chi)|> s}{(x,\chi)\in G\times\widehat{G}}}{ \int} dx d\chi.
		\end{equation}
		Then for $1<p \leq 2 $ and $T \in \mathcal{B}_p(L^2(G))$, we have
		$$\left( \int_{G\times\widehat{G}} |W^{-1}(T)(x,\chi)|^p \psi(x,\chi)^{2-p}dx\ d\chi\right)^\frac{1}{p} \lesssim M_\psi ^{\frac{2-p}{p}}\|T\|_{\mathcal{B}_p(L^2(G))}.$$
	\end{thm}
	\begin{proof}
		On $G\times\widehat{G}$, define a measure $\mu$ as
		$$\frac{\mu(x,\chi)}{dx d\chi}:= \psi^2(x,\chi).$$
		For $1<p\leq \infty$, we shall denote by $L^p(G\times\widehat{G},\mu)$, the usual $L^p$ space on $G\times\widehat{G}$ with respect to the measure $\mu,$ i.e., $$L^p(G\times\widehat{G}, \mu)=\left\{f:G  \times \widehat{G}  \to \mathbb{C}: \|f\|_{p,\mu}:=\left(\int_{G  \times \widehat{G}} |f|^p d \mu \right)^\frac{1}{p}< \infty \right\}.$$ For $T\in \mathcal{B}_p(L^2(G)),$ define a function $\Phi_T:G\times\widehat{G}\rightarrow\mathbb{C}$ as 
		$$\Phi(T)(x,\chi):=\Phi_T(x,\chi):=\frac{|W^{-1}(T)(x,\chi)|}{\psi(x,\chi)}.$$
		We now claim that $\Phi$ is a well defined sublinear bounded map from $\mathcal{B}_p(L^2(G))$ to $L^p(G\times\widehat{G},\mu)$ for $1 \leq p\leq 2 $. We will be using Marcinkiewicz's interpolation theorem (Theorem \ref{NCommCommMIT}) to prove this. In particular, we claim that $\Phi$ satisfies equation (\ref{eqn 19}) for $p_0=1$ and $p_1=2$  which follows exactly as in Theorem \ref{PIWT} and the following inequalities hold:
		\begin{align*}
			d_{\Phi_{T}}(\alpha) =&\mu\left\{(x,\chi) \in G\times\widehat{G}: |\Phi_T(x,\chi)| >  \alpha \right\} \leq \left( \frac{ \|T\|_{\mathcal{B}_2(L^2(G))}}{\alpha}\right)^2, \\
			&\mu\left\{(x,\chi) \in G\times\widehat{G}: |\Phi_T(x,\chi)| > \alpha \right\} \lesssim   \frac{M_\psi \|T\|_{\mathcal{B}_1(L^1(G))}}{\alpha}.
		\end{align*}
 Now by letting $\frac{1}{p}=1 - \theta +\frac{\theta}{2},\ 0<\theta<1,$ and applying Theorem \ref{NCommCommMIT}, we  get
		$$\left(\underset{G\times\widehat{G}}{\int} |\Phi_T(x,\chi)|^p\ d \mu(x,\chi)\right)^\frac{1}{p}\lesssim M_\psi^{\frac{2-p}{p} }\|T\|_{\mathcal{B}_p(L^2(G))}$$ which in turn gives
		\begin{align*}
			\underset{G\times\widehat{G}}{\int} \left(|W^{-1}(T)(x,\chi)|^p |\psi(x,\chi)|^{2-p}\ dxd\chi\right)^\frac{1}{p} \lesssim & M_\psi^{\frac{2-p}{p} }\|T\|_{\mathcal{B}_p(L^2(G))}. \qedhere
		\end{align*}
	\end{proof}
	As a consequence of the Paley inequality and the Hausdorff-Young inequality, we obtain the following Hausdorff-Young Paley inequality for the inverse Weyl transform.
	\begin{thm} \label{HYP}

 Consider a positive function $\psi \in L^{1, \infty}(G\times\widehat{G})$ and let $M_\psi$ be as in equation \eqref{Paley-Eqn}.
		Let $1<p \leq 2$,  and $1<p\leq b \leq p' \leq \infty$, then for $T\in \mathcal{B}_p(L^2(G))$, we have
		$$\left( \underset{G\times\widehat{G}}{\int} \left(|W^{-1}(T)(x,\chi)| \psi(x,\chi)|^{\frac{1}{b}-\frac{1}{p'}}\right)^b dx d\chi\right)^\frac{1}{b} \lesssim M_\psi^{\frac{1}{b}-\frac{1}{p'}} \|T\|_{\mathcal{B}_p(L^2(G))}.$$
	\end{thm}
	Here is the analogue of the H\"{o}rmander's theorem for the inverse Weyl transform. Let $\mathcal{F}(L^2(G))$ denote the space of all finite rank operators on $L^2(G).$
	\begin{thm}
	 Let $ 1<p \leq 2 \leq q < \infty $ and  $g \in L^{r,\infty}(G \times\widehat{G})$ where $\frac{1}{r}=\frac{1}{p}-\frac{1}{q}$.
  Consider the operator $\phi_g$ defined on $\mathcal{F}(L^2(G))$ as  $W^{-1}(\phi_g(T))=gW^{-1}(T)$. Then $\phi_g$ can be extended to a bounded operator from $\mathcal{B}_p(L^2(G))$ to $\mathcal{B}_q(L^2(G))$ and we have 
		\begin{equation} \label{eq1}
			\|\phi_g(T)\|_{\mathcal{B}_q(L^2(G))} \lesssim \sup_{s>0}s\left( \underset{|g(x, \chi)| > s}{\underset{(x, \chi ) \in G \times\widehat{G} }{\int}} dx d \chi \right)^{\frac{1}{p}-\frac{1}{q}} \|T\|_{\mathcal{B}_p(L^2(G))}.
		\end{equation}
	\end{thm}
	\begin{proof}
		Let us first assume that $p \leq q'.$ Then for $T \in \mathcal{B}_p(L^2(G))$,
		\begin{equation*} 
			\|\phi_g(T)\|_{\mathcal{B}_q(L^2(G))} \leq \|W^{-1}(\phi_g(T))\|_{q'}=\|gW^{-1}(T)\|_{q'}.
		\end{equation*}
		Now we will apply Hausdorff-Young Paley inequality of Theorem \ref{HYP} by taking $b=q'$ and $\psi(x, \chi )=(|g(x, \chi )|)^r$.  Since $\frac{1}{q'}-\frac{1}{p'}=\frac{1}{p}-\frac{1}{q}=\frac{1}{r},$ we obtain
		\begin{eqnarray*} 
			 \bigg( \underset{G \times\widehat{G}}{ \int}|W^{-1}(T)(x, \chi)g(x, \chi )|^{q'}dx  d \chi \bigg)^{\frac{1}{q'}} \lesssim \underset{s>0}{\sup}s\ \left( \underset{|g(x,\chi)| > s}{\underset{(x, \chi) \in G \times\widehat{G} }{\int}} dx  d \chi \right)^{\frac{1}{p}-\frac{1}{q}} \|T\|_{\mathcal{B}_p(L^2(G))}.
		\end{eqnarray*}
		Now by combining the above inequalities, we get the required inequality \eqref{eq1}.
		For the case when $q'<p$, one can work with $\phi^*$ whose corresponding associated function will be $\overline{g}$ and $|g|=|\overline{g}|$.
	\end{proof}

\section{Hardy-Littlewood inequality} \label{HLI}
    In this section, we study Hardy-Littlewood inequality, both scalar as well as the vector-valued case. 

    We shall begin this section by proving the scalar version of the Hardy-Littlewood inequality. This inequality is an application of the Paley inequality (Theorem \ref{Paley}).
	\begin{lem}
		Let $1<p \leq 2$. Assume that a positive function $(x,\chi) \mapsto  \mu_{(x,\chi)}$ on $G\times\widehat{G},$ have adequate rapid growth, that is, \begin{equation}
		    \underset{G\times\widehat{G}}{\int} \frac{dx d \chi}{|\mu_{(x,\chi)}|^\beta}   < \infty , \beta > 0. \label{Hardy}
		\end{equation}

		Then the following inequality holds true.
		$$ \left( \underset{G\times\widehat{G}}{\int}|\mu_{(x,\chi)}|^{-\beta(2-p)}|W^{-1}(T)(x,\chi)|^pdx d \chi  \right)^{\frac{1}{p}}\lesssim \|T\|_{\mathcal{B}_p(L^2(G))}.$$
	\end{lem}
	\begin{proof}
		We will be using Paley's inequality to prove this result.  We claim that the function $\psi(x, \chi)=|\mu_{(x, \chi)}|^{-\beta}$ satisfies equation \eqref{Paley-Eqn} required in Theorem \ref{Paley} . For $s>0$, it can be seen that 
		\begin{equation*}
			s\underset{|\mu_{(x, \chi)}|^{-\beta} > s}{ \int} dx d \chi\leq \underset{\frac{1}{s} > |\mu_{(x, \chi)}|^\beta}{ \int }\frac{dx d \chi }{|\mu_{(x, \chi)}|^\beta}\leq \int_{G \times \widehat{G}}\frac{dx d \chi}{|\mu_{(x, \chi)}|^\beta} .
		\end{equation*}
		By assumption,$$C:=\int_{G \times \widehat{G}}\frac{dx d \chi}{|\mu_{(x, \chi)}|^\beta} $$ is known to be finite.
		Hence $M_\psi \leq C< \infty$ and our claim is proved. Now by applying Theorem \ref{Paley} to the function $\psi(x)=|\mu_{(x, \chi)}|^{-\beta}$, we get the required estimate.
	\end{proof}
	
	\begin{defn}
		Let $G$ be a locally compact abelian group. Let $1<p \leq 2 \leq q < \infty$ and let $\omega: G \times \widehat{G} \to (0, \infty)$ be a weight. We say that a Banach space $X$ is of Weyl $HL$-$\omega$ type $p$ on $G$ if  there exists a constant $C> 0$  such that 
		$$  \|f\|_{L^p(G \times \widehat{G},\omega^{-(2-p)} ;X)}  \leq C \|W(f)\|_{\mathcal{B}_p[L^2(G);X]}.  $$ 
		Similarly, we say that a Banach space $X$ is of Weyl $HL$-$\omega$ cotype $q$ on $G$ if there exists a constant $C> 0$  such that 
		$$  \|W(f)\|_{\mathcal{B}_q[L^2(G);X]}  \leq C \|f\|_{L^q(G \times \widehat{G} ,\omega^{(q-2)}  ;X)}  .$$ 
	\end{defn}
\begin{thm}
	\label{HL}	X has Weyl  $HL$-$\omega$ type $p$ on $G$ if and only if $X^*$ has Weyl  $HL$-$\omega$ cotype $p'$ on $G$. \end{thm}
	\begin{proof}
		Let $X$ has Weyl  $HL$-$\omega$ type $p$, i.e,
		$\exists C>0$ such that 
		\begin{equation} \label{eqn 11}
		    \|f\|_{L^p(G \times \widehat{G},\omega^{-(2-p)} ;X)} \leq C  \|W(f)\|_{\mathcal{B}_p[L^2(G);X]} .
		\end{equation}
		
		To show that $X^*$ has Weyl  $HL$-$\omega$ cotype $p'$, we first show that given elementary tensor  $f \otimes x^* \in L^{p'}(G \times \widehat{G},\omega^{(p'-2)} ) \otimes X^*$, we have, $W(f) \otimes x^* \in \mathcal{B}_{p'}[L^2(G);X^*] $ so that $W(L^{p'}(G \times \widehat{G},\omega^{(p'-2)} ; X^* ) \subset \mathcal{B}_{p'}[L^2(G);X^*].$ Since $$\mathcal{B}_{p'}[L^2(G);X^*] =\left( \mathcal{B}_{p}[L^2(G);X]\right)^*,$$ let $T \in  \mathcal{B}_{p}[L^2(G);X] $. For $1 \leq i \leq n$, let $S_i \in B_p(L^2(G))$ and $x_i \in X$ such $\sum_{i=1}^{n} S_i \otimes x_i$ is one of the representations of $T$.
		Then, \begin{eqnarray*}
			& & \bigg |\langle W(f\otimes x^*), \sum_{i=1}^{n} S_i \otimes x_i \rangle \bigg | \\ &=& \bigg | \langle f\otimes x^*,\sum_{i=1}^{n} W^{-1}(S_i) \otimes x_i\rangle\bigg |\\
			&=&\bigg |  \sum_{i=1}^{n} \langle f, W^{-1}(S_i)\rangle \langle x^*,x_i\rangle \bigg | \\
			&\leq &  \|f\|_{L^{p'}(G \times \widehat{G},\omega^{(p'-2)} )} \|x^*\| \left(\sum_{i=1}^{n} \|W^{-1}(S_i)\|_{L^{p}(G \times \widehat{G},\omega^{-(2-p)} ) }\|x_i\|   \right).
		\end{eqnarray*}
		Since this is true for any arbitrary representation, taking infimum over all such representations we get 
		\begin{eqnarray} \label{eqn 7}
			\big |\langle W(f\otimes x^*), T \rangle| &  \leq & \|f\|_{L^{p'}(G \times \widehat{G},\omega^{(p'-2)} )} \|x^*\|\left( \|W^{-1}T\|_{L^{p}(G \times \widehat{G},\omega^{-(2-p)};X)}   \right)\\
			&\leq &  C   \|f \otimes x^*\|_{L^{p'}(G \times \widehat{G},\omega^{(p'-2)} ) \otimes X^*}  \|T\|_{\mathcal{B}_p[L^2(G);X]}. \nonumber
		\end{eqnarray}
		where we have used inequality (\ref{eqn 11}) in the last equation. Now, by using duality, we get the desired containment. \\
		Now, let $f \in L^{p'}(G \times \widehat{G},\omega^{(p'-2)};X^* ) .$ Let $\epsilon>0.$ By duality , there exists $T \in  \mathcal{B}_{p}[L^2(G);X] $ such that $\|T\|_{\mathcal{B}_{p}[L^2(G);X] }=1$ and $$\|
		W(f)\|_{ \mathcal{B}_{p'}[L^2(G);X^*])}\leq (1+\epsilon)|\langle W (f), T \rangle|.$$  Thus, 
		\begin{eqnarray*}
			\| W(f)\|_{ \mathcal{B}_{p'}[L^2(G);X^*])} & \leq & (1+\epsilon) | \langle W (f), T \rangle| \\
			&=&  (1+\epsilon) \left| \langle f, W^{-1}(T) \rangle \right| \\ 
			&\leq& (1+\epsilon)  \|f\|_{L^{p'}(G\times\widehat{G},\omega^{(p'-2)},X^* )}\|W^{-1}T\|_{L^{p}(G \times \widehat{G},\omega^{(p-2)};X)} \\
  &\leq&  (1+\epsilon) C\|f\|_{L^{p'}(G\times\widehat{G},\omega^{(p'-2)},X^* )}  \|T\|_{\mathcal{B}_p[L^2(G);X]} \\
			&=& (1+\epsilon) C\|f\|_{L^{p'}(G\times\widehat{G},\omega^{(p'-2)},X^* )} .
		\end{eqnarray*}
		As $\epsilon>0$ is arbitrary, letting $\epsilon\rightarrow 0,$ we get the desired inequality.
	\end{proof}

 \begin{cor}
 	Let G be an abelian group and $1<p \leq 2$. Assume that a positive function $(x, \chi) \to  \mu_{(x, \chi)}$ on $G \times \widehat{G}$ satisfies equation (\ref{Hardy}). Then 
  	$$\|T\|_{\mathcal{B}_{p'}(L^2(G))} \lesssim \left( \int_G|\mu_{(x, \chi)}|^{\beta(p'-2)}|W^{-1}(T)(x)|^{p'}dx d \chi \right)^{\frac{1}{p'}} .$$
 \end{cor}
	  \section{Paley-type Inequalities} \label{PTI}
 
If $f \in L^p(\mathbb{R}^{2n})$, then it  is a well know fact that $W(f) \in \mathcal{B}_{p'}(L^2(\mathbb{R}^n)).$ Paley in \cite{Paley} showed that there is an extension of Hausdorff-Young Theorem to Lorentz spaces as follows, i.e if $f \in L^p(\mathbb{T})$, then $\hat{f} \in l^{p^{\prime}, p}(\mathbb{Z})$. The Corollary \ref{cor} gives such an extension for the Weyl transform. The following is another version of the same result.

 \begin{thm}
    \label{PWT} If $f \in L^{p,p'}(G \times \widehat{G})$, then $W(f) \in \mathcal{B}_{p'}(L^2(G))$ and there exists $C>0 $ such that
     $$\|W(f)\|_{\mathcal{B}_{p'}(L^2(G))} \leq C \|f\|_{p,p'} .$$
 \end{thm}
 \begin{proof}
In view of the Marcinkiewicz interpolation theorem (\cite[Theorem 4.13]{BS}), define an operator $T$ from $L^{1}+L^{2,1}(G \times \widehat{G})$ to $\mathcal{M}_0(\mathbb{N})$   by 
 $T(f)=S_n(W(f))$.
  Since the Weyl transform maps $L^1(\mathbb{R}^{2n})$ to $\mathcal{B}_\infty(L^2(\mathbb{R}^n))$ and  $L^2(\mathbb{R}^{2n})$ to $\mathcal{B}_2(L^2(\mathbb{R}^n))$, the operator $T$ satisfies the assumptions of the Marcinkiewicz interpolation theorem. Hence $T$  maps $L^{p,p'}(G \times \widehat{G})$ to $l^{p'}(\mathbb{N})$  continuously and $$\|\{S_n(W(f))\}\|_{l^{p'}} \leq C \|f\|_{p,p'} .$$
 Since $\|W(f)\|_{\mathcal{B}_{p'}(L^2(\mathbb{R}^n))} =\|\{S_n((W(f)))\}\|_{l^{p'}}$, we get the desired result.
 \end{proof}

  \begin{defn}
		Let $G$ be a locally compact abelian group. Let $1<p \leq 2 \leq q < \infty$. We say that a Banach space $X$ is of Weyl-Paley type $p$ on $G$ if  there exist a constant $C>0$  such that 
		$$\|W(f)\|_{B_{p'}[L^2(G);X]} \leq C \|f\|_{L^{p,p'}(G \times \widehat{G} ;X)} .$$ 
		Similarly, we say that a Banach space $X$ is of Weyl-Paley cotype $q$ on $G$ if there exist a constant $C>0$  such that 
		$$\|f\|_{L^{q,q'}(G \times \widehat{G} ;X)} \leq C  \|W(f)\|_{B_{q'}[L^2(G);X]} .$$ 
	\end{defn}
\begin{thm}
    A Banach space $X$ has Weyl-Paley type $p$  on $G$ if and only if $X^*$ has Weyl-Paley cotype $p'$ on $G$.
\end{thm}
 \begin{proof}
    We simply show that if $X$ has Weyl-Paley type $p$ on $G$, then for $W(f) \otimes x^* \in \mathcal{B}_{p}(L^2(G)) \otimes X^*,$ we have  $ f \otimes x^* \in L^{p',p}(G \times \widehat{G}, X^* )$  and the rest follows as in  Theorem \ref{HL}.

   Let $u \in L^{p,p'}(G \times \widehat{G}, X)$ and let $g_i \in  L^{p,p'}(G \times \widehat{G})$ and $x_i \in X, 1 \leq i \leq n$  be such that $u=\sum_{i=1}^{n} g_i \otimes x_i$  is one of the representation of $u$.
		Then, \begin{eqnarray*}
			& & \big |\langle W^{-1}(W(f)\otimes x^*), \sum_{i=1}^{n} g_i \otimes x_i \rangle| \\ &=& \big | \langle W(f)\otimes x^*,\sum_{i=1}^{n} W(g_i) \otimes x_i\rangle\big |\\
			&=&\big |  \sum_{i=1}^{n}  \langle W(f), W(g_i)\rangle  \langle x^*,x_i \rangle \big | \\
			&\leq &  \|W(f)\|_{\mathcal{B}_{p}(L^2(G)} \|x^*\| \left(\sum_{i=1}^{n} \|W(g_i)\|_{\mathcal{B}_{p'}(L^2(G)}\|x_i\|   \right).
		\end{eqnarray*}
		Since this is true for any arbitrary representation, taking infimum over all such representations we get 
		\begin{eqnarray*}
			\big |\langle W^{-1}(W(f)\otimes x^*), u \rangle| &  \leq & \|W(f)\|_{\mathcal{B}_{p}(L^2(G)} \|x^*\| \left( \|Wu\|_{{\mathcal{B}_{p'}[L^2(G);X]}}   \right)\\
   &=& \|W(f) \otimes x^*\|_{\mathcal{B}_{p}[L^2(G);X^*]}  \left( \|Wu\|_{{\mathcal{B}_{p'}[L^2(G);X]}}   \right)\\
			&\leq &  C   \|W(f) \otimes x^*\|_{\mathcal{B}_{p}[L^2(G);X^*]}  \|u\|_{L^{p,p'}(G \times \widehat{G} ;X)}
		\end{eqnarray*}
  where we have used $X$ has Weyl-Paley type $p$ on $G$ in the last step.
\end{proof}

\begin{cor}
  Let $1 \leq p \leq 2$. Let $f$ be such that $W(f) \in \mathcal{B}_{p}(L^2(G))$, then $f \in L^{p',p}(G \times \widehat{G}) $ and there exists $C>0$ such that  
  $$\|f\|_{ L^{p',p}(G \times \widehat{G})} \leq C \|W(f)\|_{\mathcal{B}_{p}(L^2(G))} .$$
\end{cor}
\begin{prop}Let $1<p<2<q< \infty$. \label{Phl}
\begin{enumerate}
    \item[(i)] If $X$ has Weyl-Paley type $p$ on $G$, then $X$ has Weyl type $p$ on $G$. 
    \item[(ii)]  If $X$ has Weyl-Paley cotype $q$ on $G$, then $X$ has Weyl cotype $q$ on $G$.
    \item[(iii)] If $X$ has Weyl type $p$ on $G$, then it has Weyl-Paley type $p_0$ for any $p_0 \in (1, p)$ and Weyl-Paley cotype $q_0$ for any $q_0 \in (p', \infty)$  on $G$.
\end{enumerate}  
\begin{proof}
    The (i) and (ii) part of the result directly follows from the fact that $L^{p',p}(G \times \widehat{G};X ) \hookrightarrow L^{p'}(G \times \widehat{G};X )  $ and $L^{q'}(G \times \widehat{G};X ) \hookrightarrow L^{q',q}(G \times \widehat{G};X )$ are continuous embeddings.
    For the last part, if $X$ has Weyl type $p$, then 
    $$W: L^{p}(G \times \widehat{G};X) \to \mathcal{B}_{p'}[L^2(G);X] $$ is bounded map. 
    For $1<p_0< p$, there exists $\theta \in (0,1)$ such that $\frac{1}{p_0}=1-\theta+\frac{\theta}{p}$. Let $1 \leq r \leq \infty.$ Then by interpolation with parameters $(\theta, r)$, we have 
    $$W: L^{p_0,r}(G \times \widehat{G};X) \to \mathcal{B}_{p_0', r}[L^2(G);X]$$ is bounded map. In particular, choose $r=p_0'$ to get the desired result. For the other part, one can work with Weyl inverse instead.
\end{proof}
\end{prop}
\begin{prop} 
Let $1<p<2<q< \infty.$ Let $\omega: \mathbb{R}^{2n} \to [0, \infty)  $ is defined as $\omega(\cdot)=|\cdot|^{-2n}.$
\begin{enumerate}
    \item[(i)] If $X$ has Weyl-Paley type $p$ on $\mathbb{R}^n$, then it has Weyl-HL-$\omega$ cotype $p'$ on $\mathbb{R}^n$. 
    \item[(ii)]  If $X$ has Weyl-Paley cotype $q$ on $\mathbb{R}^n$, then it has Weyl-HL-$\omega$ type $q'$ on $\mathbb{R}^n$.
    \item[(iii)] If $X$ has Weyl type $p$ on $\mathbb{R}^n$, then it has Weyl-HL-$\omega $ cotype $q_0$ for any $q_0 \in (p, \infty)$ and Weyl-HL-$\omega $ type $p_0$ for any $p_0 \in (1,p)$ on $\mathbb{R}^n$.
\end{enumerate}   
\end{prop}
\begin{proof} \begin{enumerate}
    \item[(i)] 
  If $X$ has Weyl-Paley type $p$, then there exists $C>0$ such that
	$$\|W(f)\|_{B_{p'}[L^2(\mathbb{R}^{n});X]} \leq C \|f\|_{L^{p,p'}(\mathbb{R}^{2n} ;X)} .$$  Since 
$$
\int_0^{\infty} f^*(t) \frac{1}{(1 / w)^*(t)} d t \leq \int_{\mathbb{R}^d} f(x) w(x) d x .
$$ by using some rearrangements properties we get
\begin{equation}
    \int_0^{\infty} t^{ p^{\prime}/p}\left(\|f(\cdot)\|_X\right)^*(t)^{p'}\frac{d t}{t} \leq C
\int_{\mathbb{R}^{2n}}\omega(\xi)^{p'-2}\|f(\xi)\|_X^{p'} d \xi   .
 \label{eqn 10}
\end{equation}
 
The left-hand side of the above equation is  nothing but $\|f\|_{L^{p,p'}(\mathbb{R}^{2n};X)}^{p'}$. 
 Hence, using the above inequalities, we get $$\|W(f)\|_{B_{p'}[L^2(\mathbb{R}^{n});X]} \leq C\|f\|_{L^{p'}\left(\mathbb{R}^{2n} ,\omega^{ (p'-2)} ; X\right)}, $$
 for all $ f \in L^{p'}\left(\mathbb{R}^{2n} ,\omega^{(p'-2)} ; X\right ) $ as required.

\item[(ii)]Now, by using duality in equation (\ref{eqn 10}), we get
$$ \|f\|_{L^{q'}\left(\mathbb{R}^{2n}  ,\omega^{ (q'-2)} ; X\right)} \leq \|f\|_{L^{q,q'}(\mathbb{R}^{2n}  ;X)}.$$ If $X$ has Weyl-cotype $q$, then using the above equation, we get 
a constant $C> 0$  such that 
		$$ \|f\|_{L^{q'}(\mathbb{R}^{2n} ,\omega^{(q'-2)} ;X)}  \leq C \|W(f)\|_{B_{q'}[L^2(\mathbb{R}^{n} );X]} .$$ 
\item[(iii)] This is a direct consequence of the above results and Proposition \ref{Phl}.
\end{enumerate}

\section*{competing interests}
    The authors declare that they have no competing interests.
  
\end{proof}
\bibliographystyle{acm}
\bibliography{refer.bib} 

\end{document}